\documentclass [12 pt] {article}

\def\half {\mbox {$\frac {1} {2} $}}

\def\dddot#1{\mathinner{\buildrel\vbox{\kern5pt\hbox{...}}\over{#1}}}
\def\ddddot#1{\mathinner{\buildrel\vbox{\kern5pt\hbox{....}}\over{#1}}}

\catcode`\@=11

\@addtoreset{equation}{section}
\catcode`\@=10

\def\d{\mbox{\rm d}}

\begin {document}

\begin {center}
{\Large \bf The Ermakov-Pinney Equation: its varied origins and the effects of the introduction of symmetry-breaking functions}\\[3 mm]
{\large RM Morris${}^{\dagger}$ \& PGL Leach${}^{\ddagger}$ \footnote{Current address: Department of Mathematics and Institute of Systems Science, Durban University of Technology, PO Box 1334, Durban 4000, Republic of South Africa}}\\[3 mm]
{\large {$\dagger$ Department of Mathematics, Durban University of Technology, PO Box 1334, Durban 4000, Republic of South Africa}\\
{$\ddagger$ Department of Mathematics and Statistics,\\ University of Cyprus, Lefkosia 1678, Cyprus}}\\[2 mm ]
Email:  \makeatletter rmcalc85@gmail.com ; leach@ucy.ac.cy \makeatother
\end {center}
 
\begin{abstract}\noindent
The Ermakov-Pinney Equation, $$\ddot{x}+\omega^2 x=\frac{h^2}{x^3},$$ has a varied provenance which we briefly delineate.  We introduce time-dependent functions in place of the $\omega^2$ and $h^2$. The former has no effect upon the algebra of the Lie point symmetries of the equation. The latter destroys the $sl(2,\Re)$ symmetry and a single symmetry persists only when there is a specific relationship between the two time-dependent functions introduced. We calculate the form of the corresponding autonomous equation for these cases.
\end{abstract}

{\bf PACS NOS:} 02.20.Sv; 02.30.Ik; o3.65.Fd

\section {Introduction}

The Ermakov-Pinney Equation,
\begin {equation}
\ddot{x}+\omega^2 x=\frac{h^2}{x^3}, \label {1.1}
\end {equation}
is usually associated with the names of Ermakov who introduced it in a paper on integrability of second-order ordinary differential equations \cite {Ermakov 80 a} and Pinney who, some 70 years later, provided the solution in one of the more succinct papers in Mathematics \cite {Pinney 50 a}.  Ermakov used (\ref {1.1}) to obtain a first integral for the equation
\begin {equation}
\ddot {y} + \omega ^ 2 y = 0 \label {1.2}
\end {equation}
by the simple expedient of taking the difference $(1.1)y - (1.2)x$, multiplying by the integrating factor $\dot {x}y - x\dot {y} $ and integrating the result to obtain the first integral\footnote {As we have commenced with autonomous equations, the result is a first integral in the technical sense.  When one considers the actual time dependence of $x (t) $, there are those who would prefer to use the word ``invariant''.}
\begin {equation}
I = \frac{1}{2}\left(\left (\dot {x} y - x\dot {y}\right) ^ 2 + \left (\frac {y} {x}\right) ^ 2\right). \label {1.3}
\end {equation}
On the other hand Pinney simply wrote the solution of (\ref {1.1}) as
\begin {equation}
x (t) = \sqrt {A u (t) ^ 2+2 Bu (t) v (t) + Cv (t) ^ 2}, \label {1.4}
\end {equation}
where the constants, $A $, $B $ and $C $, are related according to $AC - B ^ 2 = h ^ 2/W $, $u (t) $ and $v (t) $ are linearly independent solutions of (\ref {1.2}) and $W $ is their, obviously nonzero, Wronskian.

Apart from the fact that it does seem to be a case of taking in one's own washing neither writer provided a motivation. This is not to claim that there was no motivation in either case. It simply was not revealed.

Equation (\ref {1.2}) became an object of more than just mathematical interest when work began on the modelling of controlled thermonuclear fusion.  A simple model was that of a charged particle moving in a radially symmetric electromagnetic field.  Naturally the electromagnetic field was time varying and so the $\omega $ was now time dependent.  This was not a particular concern for already at the Solvay Congress of 1911 Lorentz had proposed\footnote{It would appear that one of those ``little'' accidents of history has occluded the reality of the evolution of the adiabatic invariant.   For a recent and detailed account see \cite{Sanchez-Soto 12 a}.}  an adiabatic invariant for the simple pendulum with a string length which was slowly increasing\footnote {The situation in which the length of the pendulum was slowly decreasing was not addressed until many years later and had an entirely different outcome \cite {Ross 79 a}.}. In an Hamiltonian context with
\begin {equation}
H = \half\left (p ^ 2 + \omega (t) ^ 2 q ^ 2\right) \label {1.5}
\end {equation}
Lorentz' adiabatic invariant was
\begin {equation}
I = \frac {p ^ 2 + \omega (t) ^ 2 q ^ 2} {2\omega (t)}. \label {1.6}
\end {equation}
However, the electromagnetic field was not slowly varying and so the invariant of Lorentz was of no practical importance.  In 1966, while on sabbatical in Heidelberg, Ralph Lewis applied the asymptotic theory developed a few years earlier by Kruskal \cite {Kruskal 62 a} to the Hamiltonian (\ref {1.5}) and found, somewhat to his surprise, that the first term of his expansion was the invariant \cite {Lewis 67 a}
\begin {equation}
I = \half\left (\left (\rho p - \dot {\rho}q\right) ^ 2 + \left (\frac {q } {\rho}\right) ^ 2\right), \label {1.7}
\end {equation}
where the function, $\rho (t) $, was a solution of
\begin {equation}
\ddot {\rho} + \omega (t) ^ 2\rho = \frac {1} {\rho ^ 3} \label {1.8}
\end {equation}
which, of course, is just the equation introduced by Ermakov some 86 years before\footnote {One has the distinct impression that all of these results were obtained without knowledge of the results of others.}. This discovery led to considerable application in Physics \cite{Lewis 68 a, Lewis 68 b, Lewis 69 a}.

In 1976 Elieser and Grey \cite {Eliezer 76 a} provided an interpretation of (\ref {1.8}) in terms of the radial equation of motion for a two-dimensional time-dependent linear oscillator.  The term on the right was a consequence of the conservation of angular momentum in the system
\begin {equation}
\ddot {\bf r} + \omega (t) ^ 2 {\bf r} = 0. \label {1.9}
\end {equation}
Unfortunately the interpretation did not extend to higher dimensions \cite {Gunther 77 a}.

The third-order equation of maximal symmetry \cite {Mahomed 90 a} has the general form
\begin {equation}
\dddot {w} + 2 a (t)\dot {w} + \dot {a} (t)w = 0 \label {1.91}
\end {equation}
and has the 10-element algebra $sp (5) $ \cite {Abraham-Shrauner 95 a}.  It possesses the integrating factor, $w (t) $, with integral
\begin {equation}
w\ddot {w} - \half\dot {w} ^ 2 + a (t)w ^ 2 = 2h ^ 2 \label {1.92}
\end {equation}
which, on the substitution $w = \rho ^ 2 $ is just (\ref {1.8}).  It is known \cite {Andriopoulos 05 a} that the solutions of all scalar equations of maximal symmetry can be written in terms of combinations of the two linearly independent solutions of the corresponding second-order equation of maximal symmetry. Thus the general solution of (\ref {1.92}) is
\begin {equation}
w (t) = A u (t) ^ 2+2 Bu (t) v (t) + Cv (t) ^ 2, \label {1.93}
\end {equation}
where $u $ and $v $ are any two linearly independent solutions of
\begin {equation}
\ddot {z} + a (t)z = 0. \label {1.94}
\end {equation}

Thus we have two rational explanations for the origin of the Ermakov-Pinney Equation.

It is well known that the Ermakov-Pinney Equation possesses the three-element algebra of Lie point symmetries $sl (2,\Re) $.  What happens to the algebra if we generalise the equation to
\begin{equation}\label{1.10}
\ddot{x}+F^2(t)x=\frac{G(t)}{x^{3}}?
\end{equation}
Moreover what happens to the interpretations concerning the origin?  It is the content of this paper to explore these questions.

\section{Symmetry Analysis of (\ref{1.10})}

We examine (\ref{1.10}) for Lie point symmetries using the add-on package SYM \cite{Andriopoulos 09 a, Dimas 05 a, Dimas 06 a, Dimas 08 a, Dimas 09 a}. 

We assume a Lie point symmetry of the form
\begin{equation}\label{2.1}
\Gamma = \tau(t,x)\partial_{t}+\xi(t,x)\partial_{x}.
\end{equation}

Firstly we quote the result for the autonomous case.  The Lie point symmetries are
\begin{eqnarray*}
&& \Gamma_1 = \partial_t  \\
&& \Gamma_2 = \sin 2Ft \partial_t + Fx \cos 2Ft \partial_x \quad \mbox {\rm and} \\
&& \Gamma_3 = \cos 2Ft\partial_t - Fx \sin 2Ft \partial_x
\end{eqnarray*}
with the Lie Brackets,
\[
 \left [\Gamma_1,\,\Gamma_2\right]_{LB} = 2F\Gamma_3,\quad \left [\Gamma_1,\,\Gamma_3\right]_{LB} = -2F\Gamma_2 \quad \mbox {\rm and} \quad \left [\Gamma_2,\,\Gamma_3\right]_{LB} = - 2F\Gamma_1
\]
which is isomorphic to the algebra $A_{3, 8} $ in the Mubarakzyanov Classification Scheme \cite { Morozov 58 a, Mubarakzyanov 63 a, Mubarakzyanov 63 b, Mubarakzyanov 63 c} and more commonly known as $sl (2,\,\Re) $.  The representation above is that of $so(2,\,1)$, the algebra of noncompact rotations in three dimensions.  Note that the value of the constant parameter, $G$, does not enter into the algebra as it is a rescalable constant. 

In the case in which $F $ is nonautonomous we have a similar set of results.  Because of the so-far arbitrary time dependence of $F (t) $  it is necessary to use SYM in interactive mode.  The general symmetry has the structure
\[
\Gamma = a (t)\partial_t + \half\dot {a} (t) x\partial_x,
\]
where $a (t) $ is a solution of the third-order equation
\[
\dddot {a} + 4F (t) ^ 2\dot { a} + 4\dot {F} (t)F (t) a = 0
\]
which has the same structure as (\ref {1.91}) and so we know that the solution is given by
\[
a (t) = A_1u (t) ^ 2 + A_2u (t)v (t) + A_3v (t) ^ 2,
\]
where $u $ and $v $ are two linearly independent solutions of
\[
\ddot {r} (t) + F (t) ^ 2  r (t) = 0.
\] 
Consequently we obtain the three  Lie point symmetries
\begin{eqnarray*}
&& \Gamma_1 = u ^ 2\partial_t + u\dot {u} x\partial_x \\
&& \Gamma_2 = uv\partial_t + \half (\dot {u}v +u\dot {v}) x\partial_x \quad \mbox {\rm and} \\
&& \Gamma_3 = v ^ 2\partial_t + v\dot {v} x\partial_x
\end{eqnarray*}
with the Lie Brackets,
\[
 \left [\Gamma_1,\,\Gamma_2\right]_{LB} = W\{u,v\}\Gamma_1,\quad \left [\Gamma_1,\,\Gamma_3\right]_{LB} = 2W\Gamma_2 \quad \mbox {\rm and} \quad \left [\Gamma_2,\,\Gamma_3\right]_{LB} = W\Gamma_3,
\]
where $W = u\dot{v} - \dot{u}v$ is the Wronskian of $u$ and $v$.
The algebra is unchanged at $A_{3, 8} $, which is to say $sl (2,\Re) $.  The presence of time dependence in $F (t) $ makes no difference to the algebra and there is no restriction upon the expression for the function.

In the fully nonautonomous case of (\ref {1.10}) the analysis is more delicate and we present it in detail.  As this is a second-order equation, there are four determining equations for the coefficient functions in the posited symmetry.  The first two of the determining equations corresponding to the coefficients of $\dot{x}^3$ and $\dot{x}^2$ relate to the geometry of the underlying space \cite{Tsamparlis 15 a} and we have that 
\begin{equation}\label{2.2}
\tau=a(t)+b(t)x 
\end{equation}
and
\begin{equation}\label{2.3}
\xi=\dot{b}x^2+c(t)x+d(t).
\end{equation}

When the constraints of the actual differential equation, (\ref{1.10}), are imposed, that is the third and fourth of the determining equations, we find that 
\begin{equation}\label{2.4}
b=d=0
\end{equation}
and 
\begin{equation}\label{2.5}
c=C_{0}+\frac{1}{2}\dot{a}.
\end{equation}

In addition we obtain two equations relating the so far unspecified functions $F$ and $G$ of (\ref{1.10}) and the function $a(t)$ which is the remaining function of the coefficient functions of the general symmetry (\ref{2.1}). These are
\begin{equation}\label{2.6}
4C_{0}G-a\dot{G}=0
\end{equation}
and
\begin{equation}\label{2.7}
2F\dot{a}+a\dot{F}+\frac{1}{2}\dddot{a}=0.
\end{equation}

From (\ref{2.6}) $a$ is specified precisely by $G$ as
\begin{equation}\label{2.8}
a=\frac{4C_{0}G}{\dot{G}}.
\end{equation}

\textbf{Remark:} This result does not have a parallel in the case of autonomous $G$ as then (\ref{2.6}) says nothing about $a$ and makes $C_{0}=0$. \\

Now that $a$ is specified it is apparent that (\ref{2.7}) is a first-order ODE for $F$. The solution is given by
\begin{equation}\label{2.9}
F=\frac{M}{a^{2}}-\frac{1}{2}\left[\frac{\ddot{a}}{a}-\frac{1}{2}\left(\frac{\dot{a}}{a}\right)^{2} \right],
\end{equation}
where $M$ is a constant of integration.

Consequently for a given pair of functions $F$ and $G$ are related according to (\ref{2.8}) and (\ref{2.9}) there exists a single symmetry given by
\begin{equation}\label{2.10}
\Gamma_{s}=\frac{4G}{\dot{G}}\partial_{t}+x\left(3-\frac{2G\ddot{G}}{\dot{G}^2}\right)\partial_{x}
\end{equation}
and it is this symmetry whereby we can seek to find a possible generalisation of the Ermakov-Lewis invariant.  The reason for the single symmetry is that the parameters in the functions $F $ and $G $ occur in the coefficient functions of the symmetry. This constriction has been noted before \cite {Leach 81 a}.

\section {The autonomous form of (\ref {1.10}) given the symmetry (\ref {2.10})}

When the functions $F $ and $G $ related according to (\ref {2.8}) and (\ref {2.9}), equation (\ref {1.10}) can be rendered autonomous by finding the transformation which converts the symmetry (\ref {2.10}) to the form $\tilde {\Gamma}_s = \partial_T $. This requires the determination of the transformation
\begin {equation}
T = f (t,\,x) \quad \mbox {\rm and} \quad X = g (t,\,x), \label {4.1}
\end {equation}
where the functions of the transformation of the solutions of
\begin {equation}
\tau\frac {\partial f} {\partial t} + \xi\frac {\partial f} {\partial x} = 1 \quad \mbox {\rm and} \quad \tau\frac {\partial g} {\partial t} + \xi\frac {\partial g} {\partial x} = 0 \label {4.2}
\end {equation}
and from (\ref {2.10}) we have
\[
\tau = \frac{4G}{\dot{G}} \quad \mbox {\rm and} \quad \xi =  x\left(3-\frac{2G\ddot{G}}{\dot{G}^2}\right).
\] 
The characteristic for (\ref {4.2}b) is $x\dot {G} (t) ^ {\half}/G (t) ^ {\mbox {$\frac {3} {4} $}} $ and the second characteristic obtained for (\ref {4.2}a) is $\mbox {$\frac {3} {4} $}\log G (t) $.  As the new variables we use
\begin {equation}
T = \mbox {$\frac {3} {4} $}\log G (t) \quad \mbox {\rm and} \quad X = x \frac {\dot {G} (t) ^ {\half}} {G (t) ^ {\mbox {$\frac {3} {4} $}}}. \label {4.3}
\end {equation}
After we calculate the derivatives in terms of the new variables, make use of (\ref {1.10}) and use the relationship between $F $ and $G $ we obtain the autonomous equation
\begin {equation}
\ddot {X} (T) + 2\dot {X} (T) + \left (1+ \frac {M} {C_0 ^ 2}\right)X (T) = \frac {16} {X (T) ^ 3} \label {4.4}
\end {equation}
which is rather reminiscent of the autonomous Ermakov-Pinney Equation except for the term in the first derivative.  By construction (\ref {4.4}) has the single Lie point symmetry, $\partial_T $.

The invariants of the first extension of the symmetry $\partial_T $ are $X $ and $\dot {X} $ (obviously the overdot now refers to differentiation to $T $), say $u $ and $v $. The reduced equation is
\begin {equation}
vv' + 2v + \left (1+ \frac {M} {C_0 ^ 2}\right) = \frac {16} {u} \label {4.5}
\end {equation}
which is an Abel's Equation of the second kind for which see Kamke \cite{Kamke 83 a}[\S 4.11]. As is the case with most Abel's equations, (\ref {4.5}) is not integrable in closed form.

\section{Discussion}

In our presentation of the Lie point symmetries of the case that both $F (t) $ and $G (t) $ are constants we obtained the symmetries
\begin{eqnarray*}
&& \Gamma_1 = \partial_t  \\
&& \Gamma_2 = \sin 2Ft \partial_t + Fx \cos 2Ft \partial_x \\
&& \Gamma_3 = \cos 2Ft\partial_t - Fx \sin 2Ft \partial_x
\end{eqnarray*}
with the Lie Brackets,
\[
 \left [\Gamma_1,\,\Gamma_2\right]_{LB} = 2F\Gamma_3,\quad \left [\Gamma_1,\,\Gamma_3\right]_{LB} = -2F\Gamma_2 \quad \mbox {\rm and} \quad \left [\Gamma_2,\,\Gamma_3\right]_{LB} = - 2F\Gamma_3
\]
by a direct application of SYM.  Had we used the relationship between the solutions of the third-order equation of maximal symmetry and its corresponding second-order equation, in this case
\[
\ddot {z} + F ^ 2z = 0
\]
for which the solution set is $\{\sin Ft,\,\cos Ft\} $, we would have obtained the symmetries
\begin {eqnarray*}
& & \Gamma_1 = \sin ^ 2Ft\partial_t + F\sin Ft\cos Ftx\partial_x \\
& & \Gamma_2 = \sin Ft\cos Ft\partial_t + \half F\left (\cos ^ 2Ft - \sin ^ 2Ft\right)x\partial_x \quad \mbox {\rm and} \\
& & \Gamma_3 = \cos ^ 2Ft\partial_t - F\sin Ft\cos Ftx\partial_x
\end {eqnarray*}
and the connection between the two representations becomes quite transparent.

\section{Conclusion}

We have seen that the introduction of partial time dependence into the autonomous Ermakov-Pinney Equation, that is the usual $\omega $ becomes $\omega (t) $, makes no difference to the algebraic properties of the equation. In principle the equation is integrable, although in practice it does rely upon one's ability to integrate the equation $\ddot {x} + \omega (t) ^ 2x = 0$. The latter has been a subject of investigation for some centuries, the results of such investigation being evident in a compendium such as that of Kamke \cite {Kamke 83 a}. It is interesting to note that in the original paper of Lewis \cite {Lewis 67 a} only a few instances of the integrability of (\ref {1.8}) were provided which does suggest that the connection with the third-order equation of maximal symmetry and its standard mode of expression of solution was not generally known in those days.

When both parameters of the Ermakov-Pinney Equation depend upon time explicitly, the situation changes dramatically. The algebra $sl (2,\Re) $ is lost as is the generality of the dependence upon time in the coefficient of $x $. There is an explicit relationship in the time dependence of the two parametric functions and it is only for this specific relationship, albeit allowing for an infinite number of functions, that there exists the possibility of a single Lie point symmetry. This symmetry is sufficient for a single reduction of order to an equation which is notorious for its lack of integrability in general.

The algebra $sl (2,\Re) $ is not generally regarded as a suitable algebra for the reduction of order due to its unfortunate Lie Brackets. As is well known, reduction of order gives better results when one uses the normal subgroup. If this not be done, a local symmetry becomes a nonlocal symmetry in the reduced equation. Further reduction using a nonlocal symmetry is generally not possible in closed form. Yet the Ermakov-Pinney Equation is explicitly integrable. How can this be the case?  The answer is simple.  On a further reduction of order the nonlocal symmetry remaining is still nonlocal, but it becomes exponential nonlocal and so the calculation for the further reduction is straightforward due to the cancellation of the exponential terms.

In the case that both parameters in the Ermakov-Pinney Equation depends explicitly upon time the advantage of explicit integrability is lost.  It is reasonable to ask whether such an equation makes any physical sense. In this respect refer to the interpretation in terms of angular momentum by Eliezer and Grey \cite {Eliezer 76 a}. Their interpretation assumed a radial force for which the angular momentum is obviously a constant. However, should the angular component of the equation of motion be of the form
\[
r\ddot {\theta} + 2\dot {r}\dot {\theta} = \frac {k (t)} {r},
\]
a formal integral would exist, namely
\[
r ^ 2\dot {\theta} - \int k (t)\d t = constant
\]
and this would lead to (\ref {1.10}).

\section*{Acknowledgements}

RMM thanks the National Research Foundation of the Republic of South Africa for the granting of a postdoctoral fellowship with grant number 93183. RMM and PGLL thank the University of Cyprus for its kind hospitality while this work was being undertaken.

\end{document}